\documentclass[11pt]{article}\include{ddefs}
\usepackage{graphics,pifont,array}
\usepackage{graphicx}
\begin{document}

\begin{titlepage}

\begin{center}

{\Large{        \mbox{   }                          \\
                \mbox{   }                          \\
                \mbox{   }                          \\
                \mbox{   }                          \\
                \mbox{   }                          \\
                \mbox{   }                          \\
                \mbox{   }                          \\
                \mbox{   }                          \\
                \mbox{   }                          \\
                \mbox{   }                          \\
                \mbox{   }                          \\
  {\textbf{THE SPECIAL CASE ${\rm{{\bf{I}}}}_{{\rm{{\bf{3}}}}}$ OF THE KHOLODENKO-SILAGADZE                \\
   MULTIPLE INTEGRAL CONSIDERED ANEW   } }          \\
               \mbox{    }                          \\
              J. A. Grzesik                         \\
           Allwave Corporation                      \\
        3860 Del Amo Boulevard                      \\
                Suite 404                           \\
           Torrance, CA 90503                       \\
                \mbox{    }                         \\
           (310) 793-9620 ext 104                   \\ 
            jan.grzesik@hotmail.com                 \\
              jan@allwavecorp.com                   \\
               \mbox{     }                         \\
               \mbox{     }                         \\
          October 24, 2014                               }    }

\end{center}

\end{titlepage}

\setcounter{page}{2}

\pagenumbering{roman}
\setcounter{page}{2}
\vspace*{+2.725in}

\begin{abstract}

\parindent=0.245in

       The nested Kholodenko-Silagadze quadrature
\begin{eqnarray}
I_{n} & = & \int_{-\infty}^{\;\infty}ds_{1}\int_{-\infty}^{\;s_{1}}ds_{2}\int_{-\infty}^{\;s_{2}}ds_{3}\cdots
\int_{-\infty}^{\;s_{2n-3}}ds_{2n-2}\int_{-\infty}^{\;s_{2n-2}}ds_{2n-1}\int_{-\infty}^{\;s_{2n-1}}ds_{2n}\times                       \nonumber   \\
        &    &  \rule{1.5cm}{0mm} \times \cos(s_{1}^{2}-s_{2}^{2})\cos(s_{3}^{2}-s_{4}^{2})\cos(s_{5}^{2}-s_{6}^{2})\cdots \times      \nonumber   \\
        &    &  \rule{1.5cm}{0mm} \times \cos(s_{2n-5}^{2}-s_{2n-4}^{2})\cos(s_{2n-3}^{2}-s_{2n-2}^{2})\cos(s_{2n-1}^{2}-s_{2n}^{2})   \nonumber   \\
      & = & \frac{2}{n!}\left(\frac{\pi}{4}\right)^{n} \;,                                                                             \nonumber
\end{eqnarray}
obtained for all integers $n\geq 1$ by an elegant but indirect argument in [\textbf{1}], is tackled anew from a uniform quadrature reduction viewpoint.
Along the way, at its first instance of real difficulty when $n=3,$ the recondite quadrature
\begin{eqnarray}
\int_{\,0}^{\;\infty}\frac{\cos(u)}{u}du\int_{\,0}^{\,u}\frac{\sin^{2}(v)}{v}dv + \int_{\,0}^{\;\infty}\frac{\sin(u)}{u}du\int_{\,0}^{\,u}\frac{\sin(v)\cos(v)}{v}dv & = &
\frac{\pi^{2}}{12}\;,   \nonumber
\end{eqnarray}
heretofore presumably unknown, receives an indirect resolution with value $\pi^{2}/12$.

\end{abstract}

\pagestyle{plain}

\parindent=0.5in

\newpage

\pagenumbering{arabic}

\pagestyle{myheadings}

\setlength{\parindent}{0pt}

\pagestyle{plain}

\parindent=0.5in

\newpage
\mbox{   }

\pagestyle{myheadings}

\markright{J. A. Grzesik    \\ the kholodenko-silagadze multiple integral $i_3$}

\section{Overview}
\vspace{-3mm}

     In [{\textbf{1}}], Kholodenko and Silagadze ({\textbf{K\&S}}), in one deft, truly virtuosic stroke of great elegance established the
compact evaluation
\begin{eqnarray}
I_{n} & = & \int_{-\infty}^{\;\infty}ds_{1}\int_{-\infty}^{\;s_{1}}ds_{2}\int_{-\infty}^{\;s_{2}}ds_{3}\cdots
\int_{-\infty}^{\;s_{2n-3}}ds_{2n-2}\int_{-\infty}^{\;s_{2n-2}}ds_{2n-1}\int_{-\infty}^{\;s_{2n-1}}ds_{2n}\times                       \nonumber   \\
        &    &  \rule{1.5cm}{0mm} \times \cos(s_{1}^{2}-s_{2}^{2})\cos(s_{3}^{2}-s_{4}^{2})\cos(s_{5}^{2}-s_{6}^{2})\cdots \times      \nonumber  \\
        &    &  \rule{1.5cm}{0mm} \times \cos(s_{2n-5}^{2}-s_{2n-4}^{2})\cos(s_{2n-3}^{2}-s_{2n-2}^{2})\cos(s_{2n-1}^{2}-s_{2n}^{2})               \\
      & = & \frac{2}{n!}\left(\frac{\pi}{4}\right)^{n} \;,                                                                             \nonumber
\end{eqnarray}
holding good for all integers $n\geq 1.$  The genesis of (1) is to be found in a series development for a certain end-point altitude associated with a
sphere which has rolled without slippage along a Cornu spiral [{\textbf{1}}, Eq. (9)].  In particular, underlying (1) has been a solution to the
indicated kinematics gotten in the form of a matrix exponential, ordered in accordance with a sinistral ascent in the Cornu spiral arc length s,
suitably normalized.
\vspace{-1mm}

     At that point, instead of forcing the issue, {\textbf{K\&S}} sidestepped a head-on attempt to evaluate $I_{n}$ by passing instead to a
spinor formulation of that same movement [{\textbf{1}}, Eqs. (10)-(23)].
End-point computation under this latter viewpoint was likewise dependent upon a matrix exponential, sinistrally ordered as before with respect
to spiral arc length s, but one which admitted exact,
closed form evaluations leading to (1).  Not entirely content to rest upon the laurels of their global success (1), {\textbf{K\&S}} pondered
over the feasibility of directly evaluating (1), and proceeded to do so in connection with the three lowest tiers $I_{1},$ $I_{2},$ and $I_{3}.$
In this progression, $I_{1}$ proved to be a virtual triviality [{\textbf{1}}, Eqs. (1)-(2)], $I_{2}$ already elicited detailed attention
[{\textbf{1}}, Eqs. (24)-(30)], while finally $I_{3}$ seemed to tax to the limit the entire {\textbf{K\&S}} analytic resource
[{\textbf{1}}, Eqs. (34)-(41)].
\vspace{-1mm}

     Indeed, such direct attacks upon low level instances of $I_{n}$ were of sufficient substance to assure that, following close upon the heels
of [{\textbf{1}}], there should appear as an American Mathematical Monthly problem [{\textbf{2}}] the first
nontrivial example $I_{2}$ of (1).  Correct solutions stating that $I_{2}=\pi^{2}/16$ duly arrived from a handful of enthusiasts,
among whose names was to be found that of the undersigned.  The solution procedure ultimately printed, due to Richard Stong, was as
thoroughly {\textit{durchkomponiert}} as it was ingenious.
\vspace{-5mm}

     And there the matter rested until a chance revisit rekindled in the undersigned a hope that extension of his own solution method, which reduces
the number of integration variables from $2n$ down to just $n-1,$ could
perhaps disclose the evaluation $I_{n}=2(\pi/4)^{n}/n!$ at all levels.  Alas, the resultant formula (9) below, while a step in the right direction, is still
very far removed from this closed form to qualify.  And, even at level $n=3,$ our method impales itself upon the need to verify {\textit{a priori}} the validity of
\begin{equation}
\int_{\,0}^{\;\infty}\frac{\cos(u)}{u}du\int_{\,0}^{\,u}\frac{\sin^{2}(v)}{v}dv + \int_{\,0}^{\;\infty}\frac{\sin(u)}{u}du\int_{\,0}^{\,u}\frac{\sin(v)\cos(v)}{v}dv  = 
\frac{\pi^{2}}{12}\;.  
\end{equation}        
In fact, (2) itself could perhaps merit consideration as an open problem.  Indeed, as is sketched in Section 5 below, the manifest logarithmic divergence
as $u\uparrow\infty$ of its first inner quadrature
\begin{equation}
\int_{\,0}^{\,u}\frac{\sin^{2}(v)}{v}dv=\frac{1}{2}\int_{\,0}^{\,u}\frac{1-\cos(2v)}{v}dv
\end{equation}
\newpage
\mbox
\newline
\newline
\newline
\newline
on the left is a matter of some modest delicacy.  At the moment of composition, while we ourselves have no direct proof of (2),
numerical experiments provide strong attestations of its veracity, and this despite its ultimate, albeit refractory convergence.
On the other hand, viewed in reverse, the emergence of (2) under the requirement that (9) duly produce $I_{3}=\pi^{3}/192$ may
itself be regarded as an indirect proof.  
\vspace{-5mm}
\section{Quadrature reduction}
\vspace{-4mm}
      With use of the characteristic function $\,U_{+}(\,x\,)\,$ of the positive reals, $\epsilon>0,$
\begin{eqnarray}
U_{+}(\,x\,) &  =  &  \frac{1}{2\pi i}\int_{\,-\infty}^{\,+\infty}\frac{\,e^{i\omega x}\,}{\,\omega-i\epsilon\,}d\omega  \nonumber  \\
             &  =  & \left\{\begin{array}{ccl}
                             1 & ; &  x>0  \\
           \frac{\,1\,}{\,2\,} & ; & x=0   \\
                             0 & ; &  x<0
                          \end{array}
                   \right.   
\end{eqnarray}
the stated integral $\,I_{n}\,$ can be written without overt attention to quadrature
range curtailment as
\begin{eqnarray}
I_{n} & = & \int_{-\infty}^{\;\infty}U_{+}(\,s_{1}-s_{2}\,)\,U_{+}(\,s_{2}-s_{3}\,)\cdots U_{+}(\,s_{2n-2}-s_{2n-1}\,)\,U_{+}(\,s_{2n-1}-s_{2n}\,)\,\times    \\
  &  & \times\,\cos(s_{1}^{2}-s_{2}^{2})\cos(s_{3}^{2}-s_{4}^{2})\cdots \cos(s_{2n-3}^{2}-s_{2n-2}^{2})\cos(s_{2n-1}^{2}-s_{2n}^{2})ds_{1}ds_{2}\cdots ds_{2n-1}ds_{2n}\;\,.
                              \nonumber  
\end{eqnarray}
And then, when all cosines have been unraveled in terms of their constituent exponentials we proceed to find that
\begin{eqnarray}
2^{n}(2\pi i)^{2n-1}I_{n} & = & \rule{-4mm}{0mm}\sum_{\sigma_{1}=\pm,\cdots,\atop{\rule{-0.5mm}{0mm}\sigma_{n}=\pm}}\int_{-\infty}^{\;\infty}
\left(\,\prod_{k=1}^{2n-1}d\omega_{k}\right)\!\left(\prod_{l=1}^{2n}ds_{l}\right)
e^{i\sigma_{1}(s_{1}^{2}-s_{2}^{2})}\frac{e^{i\omega_{1}(s_{1}-s_{2})}}{\omega_{1}-i\epsilon}\frac{e^{i\omega_{2}(s_{2}-s_{3})}}{\omega_{2}-i\epsilon}e^{i\sigma_{2}(s_{3}^{2}-s_{4}^{2})} \times     \nonumber    \\
      &   &\rule{-5mm}{0mm}\frac{e^{i\omega_{3}(s_{3}-s_{4})}}{\omega_{3}-i\epsilon} \frac{e^{i\omega_{4}(s_{4}-s_{5})}}{\omega_{4}-i\epsilon}
                    \cdots e^{i\sigma_{n-2}(s_{2n-5}^{2}-s_{2n-4}^{2})}\frac{e^{i\omega_{2n-5}(s_{2n-5}-s_{2n-4})}}{\omega_{2n-5}-i\epsilon}
                   \frac{e^{i\omega_{2n-4}(s_{2n-4}-s_{2n-3})}}{\omega_{2n-4}-i\epsilon}   \times     \\
    &     &\rule{-5mm}{0mm}  e^{i\sigma_{n-1}(s_{2n-3}^{2}-s_{2n-2}^{2})}\frac{e^{i\omega_{2n-3}(s_{2n-3}-s_{2n-2})}}{\omega_{2n-3}-i\epsilon}
\frac{e^{i\omega_{2n-2}(s_{2n-2}-s_{2n-1})}}{\omega_{2n-2}-i\epsilon}
e^{i\sigma_{n}(s_{2n-1}^{2}-s_{2n}^{2})}\frac{e^{i\omega_{2n-1}(s_{2n-1}-s_{2n})}}{\omega_{2n-1}-i\epsilon}    \;.      \nonumber   
\end{eqnarray}
Further progress is made possible by square completion in the several exponents, as suggested by
\begin{eqnarray}
\sigma_{1}(s_{1}^{2}-s_{2}^{2})+\omega_{1}(s_{1}-s_{2})+\omega_{2}(s_{2}-s_{3})+\sigma_{2}(s_{3}^{2}-s_{4}^{2})+\omega_{3}(s_{3}-s_{4})+\omega_{4}(s_{4}-s_{5})+\sigma_{3}(s_{5}^{2}-s_{6}^{2})\cdots &  = &   \nonumber  \\
  &  & \rule{-12.40cm}{0mm}  \sigma_{1}(s_{1}+\sigma_{1}\omega_{1}/2)^{2}-\sigma_{1}\omega_{1}^{2}/4                                                                                   \nonumber  \\
  &  & \rule{-12.90cm}{0mm} -\sigma_{1}(s_{2}+\sigma_{1}(\omega_{1}-\omega_{2})/2)^{2}+\sigma_{1}(\omega_{1}-\omega_{2})^{2}/4                                                         \nonumber  \\ 
  &  & \rule{-12.90cm}{0mm} +\sigma_{2}(s_{3}-\sigma_{2}(\omega_{2}-\omega_{3})/2)^{2}-\sigma_{2}(\omega_{2}-\omega_{3})^{2}/4                                                         \nonumber  \\
  &  & \rule{-12.90cm}{0mm} -\sigma_{2}(s_{4}+\sigma_{2}(\omega_{3}-\omega_{4})/2)^{2}+\sigma_{2}(\omega_{3}-\omega_{4})^{2}/4                                                         \nonumber  \\
  &  & \rule{-12.90cm}{0mm} +\sigma_{3}(s_{5}-\sigma_{3}(\omega_{4}-\omega_{5})/2)^{2}-\sigma_{3}(\omega_{4}-\omega_{5})^{2}/4                                                                    \\
  &  & \rule{-12.90cm}{0mm} -\sigma_{3}(s_{6}+\sigma_{3}(\omega_{5}-\omega_{6})/2)^{2}+\sigma_{3}(\omega_{5}-\omega_{6})^{2}/4                                                         \nonumber  \\
  &  &             \rule{-10cm}{0mm}                    \cdots\cdots                                                                                                                    \nonumber  \\
  &  & \rule{-12.90cm}{0mm} +\sigma_{n}(s_{2n-1}+\sigma_{n}(\omega_{2n-1}-\omega_{2n-2})/2)^{2}-\sigma_{n}(\omega_{2n-1}-\omega_{2n-2})^{2}/4                                          \nonumber  \\
  &  & \rule{-12.90cm}{0mm} -\sigma_{n}(s_{2n}+\sigma_{n}\omega_{2n-1}/2)^{2}+\sigma_{n}\omega_{2n-1}^{2}/4\;.                                                                         \nonumber
\end{eqnarray}
It becomes clear that all $2n$ quadratures over variables $s_{l}$ can be performed at once in closed form whereas, insofar as variables
$\omega_{k}$ are concerned, their squares appear only at even-indexed spots $\omega_{2m}.$  So
\newpage
\mbox{    }
\newline
\begin{eqnarray}
2^{n}(2\pi i)^{2n-1}I_{n} & = & \pi^{n}\!\!\!\sum_{\sigma_{1}=\pm,\cdots,\atop{\rule{-0.5mm}{0mm}\sigma_{n}=\pm}}\int_{-\infty}^{\;\infty}\!
\frac{e^{-i\sigma_{1}\omega_{1}\omega_{2}/2}}{\omega_{1}-i\epsilon}d\omega_{1}\!\left\{
\prod_{k=2}^{n-1}\frac{e^{-i\sigma_{k}\omega_{2k-1}(\omega_{2k}-\omega_{2k-2})/2}}{\omega_{2k-1}-i\epsilon}d\omega_{2k-1}\right\}
\frac{e^{i\sigma_{n}\omega_{2n-1}\omega_{2n-2}/2}}{\omega_{2n-1}-i\epsilon}d\omega_{2n-1} \nonumber   \\
      &    & \rule{4.4cm}{0mm}\times \left\{\prod_{k=1}^{n-2}\frac{e^{i(\sigma_{k}-\sigma_{k+1})\omega_{2k}^{2}/4}}{\omega_{2k}-i\epsilon}d\omega_{2k}\right\}
\frac{e^{-i\sigma_{n}\omega_{2n-2}^{2}/4}}{\omega_{2n-2}-i\epsilon}d\omega_{2n-2}  \;,
\end{eqnarray}
whereupon reference to (4) gives
\begin{eqnarray}
2^{n}(2\pi i)^{n-1}I_{n} & = &  \pi^{n}\sum_{\sigma_{1}=\pm,\cdots,\atop{\rule{-0.5mm}{0mm}\sigma_{n}=\pm}}\int_{-\infty}^{\;\infty}\!
U_{+}(-\sigma_{1}\omega_{2})\left\{\,\prod_{k=2}^{n-1}U_{+}\!\left(\rule{0mm}{3.25mm}\sigma_{k}(\omega_{2k-2}-\omega_{2k})\right)\right\}U_{+}(\sigma_{n}\omega_{2n-2}) \times  \nonumber   \\
    &    &  \rule{4.1cm}{0mm} \times \left\{\prod_{k=1}^{n-2}\frac{e^{i(\sigma_{k}-\sigma_{k+1})\omega_{2k}^{2}/4}}{\omega_{2k}-i\epsilon}d\omega_{2k}\right\}
\frac{e^{-i\sigma_{n}\omega_{2n-2}^{2}/4}}{\omega_{2n-2}-i\epsilon}d\omega_{2n-2}.
\end{eqnarray}

\section{The special case with \mbox{\boldmath${\rm{n}}=3$}}
\vspace{-3mm}

We now have, with two integration variables $\omega_{2}$ and $\omega_{4}$ as opposed to the initial $\{s_{k}\}_{k=1}^{6},$
\begin{eqnarray}
I_{3} & = & -\frac{\pi}{32}\sum_{\sigma_{1}=\pm,\sigma_{2}=\pm,\atop{\rule{-0.2mm}{0mm}\sigma_{3}=\pm}}\int_{-\infty}^{\;\infty}\!
U_{+}(-\sigma_{1}\omega_{2})U_{+}\!\left(\rule{0mm}{3.25mm}\sigma_{2}(\omega_{2}-\omega_{4})\right)U_{+}(\sigma_{3}\omega_{4})
\frac{e^{i(\sigma_{1}-\sigma_{2})\omega_{2}^{2}/4}}{\omega_{2}-i\epsilon}d\omega_{2}
\frac{e^{-i\sigma_{3}\omega_{4}^{2}/4}}{\omega_{4}-i\epsilon}d\omega_{4}   \nonumber   \\
        & = & -\frac{\pi}{32}\left[\rule{0mm}{7mm}\int_{\,0}^{\;\infty}\!\frac{e^{-i\omega_{4}^{2}/4}}{\omega_{4}-i\epsilon}d\omega_{4}\int_{\,\omega_{4}}^{\;\infty}\!
                                             \frac{e^{-i\omega_{2}^{2}/2}}{\omega_{2}-i\epsilon}d\omega_{2} + 
               \int_{\,0}^{\;\infty}\!\frac{e^{-i\omega_{4}^{2}/4}}{\omega_{4}-i\epsilon}d\omega_{4}\int_{\,0}^{\,\omega_{4}}\!
                                             \frac{1}{\omega_{2}-i\epsilon}d\omega_{2}\;+ \right.                     \nonumber  \\
        &   & \rule{3.5cm}{0mm}   + \int_{\,0}^{\;\infty}\!\frac{e^{-i\omega_{4}^{2}/4}}{\omega_{4}-i\epsilon}d\omega_{4}\int_{\,-\infty}^{\,0}\!
                                             \frac{e^{i\omega_{2}^{2}/2}}{\omega_{2}-i\epsilon}d\omega_{2}\,+            \nonumber  \\
        &   &                                                                                                         \nonumber           \\
        &   & \rule{0.7cm}{0mm} + \int_{-\infty}^{\,0}\!\frac{e^{i\omega_{4}^{2}/4}}{\omega_{4}-i\epsilon}d\omega_{4}\int_{-\infty}^{\,\omega_{4}}\!
                                             \frac{e^{i\omega_{2}^{2}/2}}{\omega_{2}-i\epsilon}d\omega_{2} + 
               \int_{-\infty}^{\,0}\!\frac{e^{i\omega_{4}^{2}/4}}{\omega_{4}-i\epsilon}d\omega_{4}\int_{\,\omega_{4}}^{\,0}\!
                                             \frac{1}{\omega_{2}-i\epsilon}d\omega_{2}\;+                             \nonumber  \\
        &   & \rule{3.5cm}{0mm}\left.     + \int_{-\infty}^{\,0}\!\frac{e^{i\omega_{4}^{2}/4}}{\omega_{4}-i\epsilon}d\omega_{4}\int_{\,0}^{\,\infty}\!
                                             \frac{e^{-i\omega_{2}^{2}/2}}{\omega_{2}-i\epsilon}d\omega_{2}  \rule{0mm}{7mm} \,\right]                 \\
        & = & -\frac{\pi}{32}\left[\rule{0mm}{7mm}\int_{\,0}^{\;\infty}\!\frac{e^{-i\omega_{4}^{2}/4}}{\omega_{4}-i\epsilon}d\omega_{4}\int_{\,0}^{\;\infty}\!
                                             \frac{e^{-i\omega_{2}^{2}/2}}{\omega_{2}-i\epsilon}d\omega_{2} - 
               \int_{\,0}^{\;\infty}\!\frac{e^{-i\omega_{4}^{2}/4}}{\omega_{4}-i\epsilon}d\omega_{4}\int_{\,0}^{\,\omega_{4}}\!
                                             \frac{e^{-i\omega_{2}^{2}/2}-1}{\omega_{2}-i\epsilon}d\omega_{2}\;+ \right.                     \nonumber  \\
        &   & \rule{3.5cm}{0mm}   + \int_{\,0}^{\;\infty}\!\frac{e^{-i\omega_{4}^{2}/4}}{\omega_{4}-i\epsilon}d\omega_{4}\int_{\,-\infty}^{\,0}\!
                                             \frac{e^{i\omega_{2}^{2}/2}}{\omega_{2}-i\epsilon}d\omega_{2}\,+            \nonumber  \\
        &   &                                                                                                       \nonumber             \\
        &   & \rule{0.7cm}{0mm} + \int_{-\infty}^{\,0}\!\frac{e^{i\omega_{4}^{2}/4}}{\omega_{4}-i\epsilon}d\omega_{4}\int_{-\infty}^{\,0}\!
                                             \frac{e^{i\omega_{2}^{2}/2}}{\omega_{2}-i\epsilon}d\omega_{2} - 
               \int_{-\infty}^{\,0}\!\frac{e^{i\omega_{4}^{2}/4}}{\omega_{4}-i\epsilon}d\omega_{4}\int_{\,\omega_{4}}^{\,0}\!
                                             \frac{e^{i\omega_{2}^{2}/2}-1}{\omega_{2}-i\epsilon}d\omega_{2}\;+                             \nonumber  \\
        &   & \rule{3.5cm}{0mm}\left.     + \int_{-\infty}^{\,0}\!\frac{e^{i\omega_{4}^{2}/4}}{\omega_{4}-i\epsilon}d\omega_{4}\int_{\,0}^{\,\infty}\!
                                             \frac{e^{-i\omega_{2}^{2}/2}}{\omega_{2}-i\epsilon}d\omega_{2}  \rule{0mm}{7mm} \,\right]  \;.        \nonumber        
 \end{eqnarray}

\newpage
\mbox{   }
\newline
\newline
\newline
Additional consolidation then takes the form
\begin{eqnarray}
I_{3}   & = & -\frac{\pi}{32}\left[\rule{0mm}{7mm}\left\{\rule{0mm}{6mm}\int_{-\infty}^{\,0}\!\frac{e^{i\omega_{4}^{2}/4}}{\omega_{4}-i\epsilon}d\omega_{4}+
                                                                        \int_{\,0}^{\,\infty}\!\frac{e^{-i\omega_{4}^{2}/4}}{\omega_{4}-i\epsilon}d\omega_{4}\right\}\!
                                                  \left\{\rule{0mm}{6mm}\int_{-\infty}^{\,0}\!\frac{e^{i\omega_{2}^{2}/2}}{\omega_{2}-i\epsilon}d\omega_{2}+
                                                                        \int_{\,0}^{\,\infty}\!\frac{e^{-i\omega_{2}^{2}/2}}{\omega_{2}-i\epsilon}d\omega_{2}\right\} \right. \nonumber  \\    
        &   &                                                                                                       \nonumber             \\
        &   &  \left.\rule{1cm}{0mm} - \int_{-\infty}^{\,0}\!\frac{e^{i\omega_{4}^{2}/4}}{\omega_{4}-i\epsilon}d\omega_{4}\int_{\,\omega_{4}}^{\,0}\!
                                             \frac{e^{i\omega_{2}^{2}/2}-1}{\omega_{2}-i\epsilon}d\omega_{2} \rule{0mm}{7mm}
                      - \int_{\,0}^{\;\infty}\!\frac{e^{-i\omega_{4}^{2}/4}}{\omega_{4}-i\epsilon}d\omega_{4}\int_{\,0}^{\,\omega_{4}}\!
                                             \frac{e^{-i\omega_{2}^{2}/2}-1}{\omega_{2}-i\epsilon}d\omega_{2}\, \right]\;.     
\end{eqnarray}

     Now, the inherent reality of $I_{3}$ is assured by noting, first, that each of the curly brackets in (11), being built up as a difference of a quantity
and its complex conjugate, is purely imaginary, so that their product is automatically real.  By the same token, the remaining two terms on the right in (11)
are recognized to be the sum of two complex conjugates and thus to compose a real contribution once more.  In performing the required identification of real and
imaginary parts, we pass to the limit $\epsilon \downarrow 0+$ and avail ourselves of the standard shorthand
\begin{equation}
\lim_{\atop \epsilon\downarrow 0+}\frac{1}{\omega\mp i\epsilon}=\lim_{\atop \epsilon\downarrow 0+}\left\{\rule{0mm}{7mm}\frac{\omega}{\omega^{2}+\epsilon^{2}}
\pm i\frac{\epsilon}{\omega^{2}+\epsilon^{2}}\right\}=\frac{P}{\omega}\pm i\pi \delta(\omega)
\end{equation}
involving the Cauchy principal value $P$ and Dirac's delta $\delta.$  Hence
\begin{eqnarray}
\int_{-\infty}^{\,0}\!\frac{e^{i\omega_{4}^{2}/4}}{\omega_{4}-i\epsilon}d\omega_{4}+
                               \int_{\,0}^{\,\infty}\!\frac{e^{-i\omega_{4}^{2}/4}}{\omega_{4}-i\epsilon}d\omega_{4}  & = &
2\,i\Im\left(\rule{0mm}{6mm}\int_{\,0}^{\,\infty}\!\left\{\cos(\omega_{4}^{2}/4)-i\sin(\omega_{4}^{2}/4)\right\}
\left\{\rule{0mm}{5mm}\frac{P}{\omega_{4}}+i\pi\delta(\omega_{4})\right\}d\omega_{4}\right)  \nonumber   \\
     &   =  &  2\,i\left(\rule{0mm}{6mm}\int_{\,0}^{\,\infty}
\left\{\rule{0mm}{5mm}\pi\cos(\omega_{4}^{2}/4)\delta(\omega_{4})-\sin(\omega_{4}^{2}/4)\frac{P}{\omega_{4}}\right\}d\omega_{4}\right)    \\
   &   =  &  2\,i\left(\rule{0mm}{5mm}\frac{\pi}{2}-\frac{\pi}{4}\right)    \nonumber \\
   &   =  &   \frac{i\pi}{2}\;,     \nonumber
\end{eqnarray}
whereas identical reasoning gives
\begin{eqnarray}
\int_{-\infty}^{\,0}\!\frac{e^{i\omega_{2}^{2}/2}}{\omega_{2}-i\epsilon}d\omega_{2}+
                               \int_{\,0}^{\,\infty}\!\frac{e^{-i\omega_{2}^{2}/2}}{\omega_{2}-i\epsilon}d\omega_{2}  & = &  \frac{i\pi}{2}     
\end{eqnarray}
so that the product of terms in the curly brackets in the first line on the right in (11) is just $-\pi^{2}/4.$

      The remaining two terms on the right in (11) augment the amount $-\pi^{2}/4$ now found with
\begin{eqnarray}
-2\Re\left(\!\rule{0mm}{6mm}\int_{\,0}^{\;\infty}\!\frac{e^{-i\omega_{4}^{2}/4}}{\omega_{4}-i\epsilon}d\omega_{4}\int_{\,0}^{\,\omega_{4}}\!
                                             \frac{e^{-i\omega_{2}^{2}/2}-1}{\omega_{2}-i\epsilon}d\omega_{2} \right)  &   =   & 
                 -2\left(\!\rule{0mm}{6mm}\int_{\,0}^{\;\infty}
                              \left\{\rule{0mm}{5mm}\cos(\omega_{4}^{2}/4)\frac{P}{\omega_{4}}+\pi\sin(\omega_{4}^{2}/4)\delta(\omega_{4})\right\}d\omega_{4}\,\times \right.\nonumber  \\
 &   & \rule{0.4cm}{0mm} \times \int_{\,0}^{\,\omega_{4}}
       \left\{\rule{0mm}{5mm}\left(\rule{0mm}{3mm}\cos(\omega_{2}^{2}/2)-1\right)\frac{P}{\omega_{2}}+\pi\sin(\omega_{2}^{2}/2)\delta(\omega_{2})\right\}d\omega_{2}  \,-   \nonumber      \\
 &   & \rule{4mm}{0mm}-\int_{\,0}^{\;\infty}
        \left\{\rule{0mm}{5mm}\pi\cos(\omega_{4}^{2}/4)\delta(\omega_{4})-\sin(\omega_{4}^{2}/4)\frac{P}{\omega_{4}}\right\}d\omega_{4}\,\times         \\
 &   & \left.\rule{0mm}{6mm}  \rule{0.4cm}{0mm} \times \int_{\,0}^{\,\omega_{4}}
       \left\{\rule{0mm}{5mm}\pi\left(\rule{0mm}{3mm}\cos(\omega_{2}^{2}/2)-1\right)\delta(\omega_{2})-\sin(\omega_{2}^{2}/2)\frac{P}{\omega_{2}}\right\}d\omega_{2}\right)\;.  \nonumber
\end{eqnarray}
\newpage
\mbox{   }
\newline
\newline
\newline                     
At this point there intervene several obvious simplifications having
\begin{eqnarray}
-2\Re\left(\!\rule{0mm}{6mm}\int_{\,0}^{\;\infty}\!\frac{e^{-i\omega_{4}^{2}/4}}{\omega_{4}-i\epsilon}d\omega_{4}\int_{\,0}^{\,\omega_{4}}\!
                                             \frac{e^{-i\omega_{2}^{2}/2}-1}{\omega_{2}-i\epsilon}d\omega_{2} \right) &  =   &
                 -2\left(\!\rule{0mm}{6mm}\int_{\,0}^{\;\infty}\frac{\cos(\omega_{4}^{2}/4)}{\omega_{4}}d\omega_{4}\int_{\,0}^{\,\omega_{4}}
                                        \frac{\cos(\omega_{2}^{2}/2)-1}{\omega_{2}}d\omega_{2}  -  \right.    \nonumber   \\
    &      & \left.  \rule{0mm}{6mm} \rule{2mm}{0mm}  - \int_{\,0}^{\;\infty}\frac{\sin(\omega_{4}^{2}/4)}{\omega_{4}}d\omega_{4}\int_{\,0}^{\,\omega_{4}}
                                        \frac{\sin(\omega_{2}^{2}/2)}{\omega_{2}}d\omega_{2} \right)            \nonumber  \\
    &   =  & \rule{2mm}{0mm}  4\left(\!\rule{0mm}{6mm}\int_{\,0}^{\;\infty}\frac{\cos(\omega_{4}^{2}/4)}{\omega_{4}}d\omega_{4}\int_{\,0}^{\,\omega_{4}}
                                        \frac{\sin^{2}(\omega_{2}^{2}/4)}{\omega_{2}}d\omega_{2}\,+  \right.       \\
    &      & \left.\rule{2mm}{0mm} +  \int_{\,0}^{\;\infty}\frac{\sin(\omega_{4}^{2}/4)}{\omega_{4}}d\omega_{4}\int_{\,0}^{\,\omega_{4}}
                                        \frac{\sin(\omega_{2}^{2}/4)\cos(\omega_{2}^{2}/4)}{\omega_{2}}d\omega_{2}    \rule{0mm}{6mm}\right)    \nonumber    \\
    &   =  & \rule{5mm}{0mm}  \int_{\,0}^{\;\infty}\frac{\cos(u)}{u}du\int_{\,0}^{\,u}
                                        \frac{\sin^{2}(v)}{v}dv   \nonumber  \\
    &      &  \rule{2.4cm}{0mm}         + \int_{\,0}^{\;\infty}\frac{\sin(u)}{u}du\int_{\,0}^{\,u}\frac{\sin(v)\cos(v)}{v}dv        \nonumber  \\
    &   =  & \rule{5mm}{0mm}  \int_{\,0}^{\;\infty}\frac{\cos(u)}{u}du\int_{\,0}^{\,u}
                                        \frac{\sin^{2}(v)}{v}dv   \nonumber  \\
    &      &  \rule{2.4cm}{0mm}         + \int_{\,0}^{\;\infty}\frac{\cos(u)\sin(u)}{u}du\int_{\,u}^{\,\infty}\frac{\sin(v)}{v}dv        \nonumber 
\end{eqnarray}
as their net outcome.  Agreement with the {\textbf{K\&S}} result is thus achieved only if we can verify that
\vspace{-1mm}
\begin{eqnarray}
\int_{\,0}^{\;\infty}\frac{\cos(u)}{u}du\int_{\,0}^{\,u}\frac{\sin^{2}(v)}{v}dv + \int_{\,0}^{\;\infty}\frac{\sin(u)}{u}du\int_{\,0}^{\,u}\frac{\sin(v)\cos(v)}{v}dv & = &
\frac{\pi^{2}}{12}\;.  \\
     &    &            \nonumber
\end{eqnarray}
Happily enough, a numerical exploration of this query via coding in {\textbf{FORTRAN}} amounted to a near confirmation, albeit not
without a certain amount of exertion in securing convergence for the first triangular integral on the left ({\textit{cf.}} Section 5 below).
\vspace{-5mm}

\section{The special cases with \mbox{\boldmath${\rm{n}}=1\,\&\,2$} noted in passing}
\vspace{-2mm}

      After the vigorous exercise chronicled in (10)-(17), the antecedent cases involving $n=1$ and $n=2$ are somewhat anticlimactic.  Needless to say,
it has been implicit all along that formula (9) applies only for $n\geq 2.$  At its lowest, $n=2$ level it now gives
\begin{eqnarray}
I_{2}&=&\frac{\pi}{8i}\sum_{\sigma_{1}=\pm,\atop{\rule{-0.7mm}{0mm}\sigma_{2}=\pm}}\int_{-\infty}^{\;\infty}\!
U_{+}(-\sigma_{1}\omega_{2})U_{+}(\sigma_{2}\omega_{2})\frac{e^{-i\sigma_{2}\omega_{2}^{2}/4}}{\omega_{2}-i\epsilon}d\omega_{2}   \nonumber \\
     &=& \frac{\pi}{8i}\left(\rule{0mm}{6mm}\int_{\,0}^{\;\infty}\frac{e^{-i\omega_{2}^{2}/4}}{\omega_{2}-i\epsilon}d\omega_{2}\,+ 
     \,\int_{-\infty}^{\,0}\frac{e^{i\omega_{2}^{2}/4}}{\omega_{2}-i\epsilon}d\omega_{2} \,\right)      \nonumber     \\
  &=& \frac{\pi}{4}\Im\left(\rule{0mm}{6mm}\int_{\,0}^{\;\infty}\left\{\cos(\omega^{2}/4)-i\sin(\omega^{2}/4)\right\}
                                               \left\{\rule{0mm}{5mm}\frac{P}{\omega}+i\pi\delta(\omega)\right\}d\omega\,\right)    \\
  & &                                                                                                                  \nonumber 
\end{eqnarray}                                                                                                             
\begin{center}
\vspace{-6mm}
        Eq. (18) continued
\vspace{-3mm}
\end{center}
\newpage
\mbox{   }
\begin{center}
        Eq. (18) continued
\vspace{-1cm}
\end{center}
\setcounter{equation}{17}
\begin{eqnarray}
     & = & \frac{\pi}{4}\left(\rule{0mm}{6mm}\frac{\pi}{2}-\int_{\,0}^{\infty}\frac{\sin(\omega^{2}/4)}{\omega}d\omega\,\right)  \\
     & = & \frac{\pi}{4}\left(\rule{0mm}{6mm}\frac{\pi}{2}-\frac{1}{2}\int_{\,0}^{\infty}\frac{\sin(\zeta)}{\zeta}d\zeta\,\right)  \nonumber  \\
     & = & \frac{\pi^{2}}{16}\;,     \nonumber
\end{eqnarray}
which is correct ({\textit{cf.} [{\textbf{1}},{\textbf{2}}]).

     And as for $I_{1},$ much of the ponderous machinery from (6)-(9) simply evaporates, leaving us with just
\begin{equation}
I_{1}  =  \frac{1}{4 i}\sum_{\sigma=\pm}\int_{-\infty}^{\;\infty}\frac{1}{\omega-i\epsilon}d\omega = \frac{\pi}{2}\;,
\end{equation}
correct once again.
\vspace{-5mm}

\section{A convergence afterthought}

     {\textit{A priori,}} the inner integral
\begin{equation}
\int_{\,0}^{\,u}\frac{\sin^{2}(v)}{v}dv=\frac{1}{2}\int_{\,0}^{\,u}\frac{1-\cos(2v)}{v}dv
\end{equation}
appearing in the first term on the left in (17), when considered in isolation, clearly diverges as $\log(u)$ when $u\uparrow\infty.$
Convergence for the iterated, composite integral is nevertheless rescued by the presence of factor $1/u,$ since, on the one hand,
L'H{\^{o}}pital's rule assures that $\log(u)/u\approx 1/u,$ whereas, for any fixed $\alpha>0,$ the convergence of
$\int_{\alpha}^{\,u}\left\{\cos(u)/u\right\}du$ as $u\uparrow\infty$ is guaranteed, if nothing else, by that of $\int_{\alpha}^{\,u}\left\{\sin(u)/u\right\}du.$
Any residual convergence anxieties around the origin, $u\downarrow 0\,\!+,$ are put to rest by inspection.

    The need to contend with $\log(u)/u$ as opposed to only $1/u$ has the heuristic effect of dilating the quadrature interval required to
secure convergence of the first, triangular integral on the left in (17).  Having enormous computing power within easy reach, we were able
to integrate up to $u_{{\rm{max}}}=2000000,$ an astounding feat performed in negligible time.  With Gauss-Legendre quadrature (GLQ)
rules at level 10 endlessly concatenated across slots of unit magnitude, we were able to produce $0.9892905538  \ldots\,  \times\, \pi^{2}/12$ on the left in (17),
whereas the slightly more coarse level 8 quadrature rules similarly gave $1.0175685832  \ldots\, \times\, \pi^{2}/12.$\footnote{Whereas
GLQ rules of too low an order are discouraged for obvious reasons, slightly more subtle considerations
intrude to oppose rule invocation in the opposite direction.  At issue is the fact that our simpleminded, plodding advance across a lattice of quadrature
nodes set down in advance has the unpleasant effect, but only as regards the inner integrals, of effectively inducing an incomplete data interpolation, lopsided toward
the lower end of each quadrature slot, the upper remainder, progressively shrinking, to be sure, being padded with values artificially set equal to zero.
Accordingly, the numerical estimates now cited are to be understood merely as {\textit{grosso modo}} indicators and no more.} 
A measure of quality control was provided by
the standard result
\begin{equation}
\int_{\,0}^{\,\infty}\frac{\sin(u)}{u}du=\frac{\pi}{2}\,,
\end{equation}
\newpage
\mbox{    }
\mbox{    }
\newline
\newline
\newline      
the numerical quadratures of which respectively gave $0.9999997745  \ldots\, \times\, \pi/2$ and $0.9999997447  \ldots\, \times\, \pi/2.$
In general, the second integral on the left in (17) emerged as dominant, being positive and having a magnitude roughly five times that of the first,
weakly convergent one, all of whose estimates were negative.  Save for the caveats implied in Footnote 1, the actual coding mechanics elicited by
the triangular forms in (17) proved to be but little more involved than those encountered during routine, one-dimensional integration.  We
resisted the temptation to explore the numerical landscape after the integrals had been transposed into a polar co\"{o}rdinate $(r,\varphi)$ setting with
$r=\sqrt{u^{2}+v^{2}}$ and $\varphi=\arctan(v/u).$

\section{A parting comment}

      One could in principle continue along this route, harnessing (9) so as to attempt an {\textit{ab initio}} calculation of
$I_{4},$ $I_{5},$ and so on {\textit{ad infinitum,}} with increasingly intricate counterparts to (17) most likely to be discovered
along the way.  But with the algebraic complexity already unleashed in connection with $I_{3}$ now in plain view, it is evident
that this Sysiphean endeavor would very quickly tend to destabilize one's mental equilibrium.  Hence it is that we can rightly
stand in awe at the success of {\textbf{K\&S}} in assembling an indirect pattern of arguments to produce their compact
evaluation $I_{n}=2(\pi/4)^{n}/n!$ at all indices $n.$  Indeed, the very regularity of that evaluation pleads for
a direct demonstration, tethered perhaps to some sort of symmetry under permutation of integration variables $s_{k}$ taken in
couplets.  Such avenues were explored to a certain extent, but without tangible success at the time of writing.

\section{References}
\parindent=0.0in
1.  A. L. Kholodenko and Z. K. Silagadze, {\textbf{When physics helps mathematics: calculation of the sophisticated multiple integral}}, arXiv:1201.1975v1 [math.CA] 10 Jan 2012.

2.  Z. K. Silagadze, Problem 11621, {\textbf{A Quadruple Integral}}, The American Mathematical Monthly, Vol. 119, No. 2, Feb. 2012, p. 161:  submitted solution by J. A. Grzesik,
March 20, 2012 (one among several contributors);  published solution by Richard Stong, {\textit{ibid.,} Vol. 120, No. 7, Aug.-Sept. 2013, pp. 667-669.

\end{document}